\documentstyle{amsart}
\newcommand{\freeprod}{{\ast_{i\in I}}G_i}
\newtheorem{theorem}{Theorem}
\newtheorem{lemma}[theorem]{Lemma}
\newtheorem{corollary}[theorem]{Corollary}
\newcounter{clause}
\newenvironment{parts}{\begin{list}{\rom{(\thetheorem.\theclause)}}%
{\usecounter{clause}\setlength{\itemindent}{0in}\setlength{\leftmargin}{0in}}}%
{\end{list}}

\begin{document}

\title[Commutators as Powers]{Commutators as Powers in Free Products of Groups}

\author[L.~P.~Comerford, Jr.]{Leo P. Comerford, Jr.}
\address{Department of Mathematics\\Eastern Illinois 
University\\Charleston, Illinois 61920\\U.~S.~A.}
\email{cflpc@@eiu.edu}

\author[C.~C.~Edmunds]{Charles C. Edmunds}
\address{Department of Mathematics\\Mount St.\ Vincent 
University\\Halifax, Nova Scotia B3M 2J6\\Canada}
\email{edmundsc@@ash.msvu.ca} 

\author[G.~Rosenberger]{Gerhard Rosenberger}
\address{Fachbereich Mathematik\\Universit\"{a}t Dortmund\\
4600 Dortmund 50\\Germany}
\email{UMA004\%DDOHRZ11.BITNET@@vm.gmd.de}

\subjclass{Primary 20E06; Secondary 20F12}

\maketitle

\begin{abstract}
The ways in which a nontrivial commutator can be a proper power in a free 
product of groups are identified.
\end{abstract}

It is well known that in a free group, a nontrivial commutator cannot be a 
proper power.  This seems to have been noted first by M.~P.~Sch\"{u}tzenberger 
\cite{Schutz59}.  It is, however, possible for a nontrivial commutator to be a 
proper power in a free product. 
Our aim in this paper is to determine the ways in which this can happen. 

\begin{theorem}\label{main}
Let $G=\freeprod$, the free product of nontrivial free factors $G_i$.  If 
$V, X, Y\in G$ and $V^m=X^{-1}Y^{-1}XY=[X,Y]$ for some $m\ge 2$, then 
either
\begin{parts}
\item\label{factor} $V\in W^{-1}G_i W$ for some $W\in G$, $i\in I$, and $V^m$ 
is a commutator in $W^{-1}G_i W$, or
\item\label{m=even} $m$ is even, $V=AB$ with $A^2=B^2=1$, and 
$V^m=[A,B(AB)^{(m-2)/2}]$, or
\item\label{m=odd} $m$ is odd, $V=AC^{-1}AC$ with $A^2=1$, and 
$V^m=[A,C(AC^{-1}AC)^{(m-1)/2}]$, or
\item\label{m=6} $m=6$, $V=AB$ with $A^2=B^3=1$, and 
$V^6=[B^{-1}ABA,B(AB)^2]$, or
\item\label{m=3} $m=3$, $V=AB$ with $A^3=B^3=1$, and $V^3=[BA^{-1},BAB]$, or
\item\label{m=2} $m=2$, $V=AB$ with $A^2=1$, and $B^{-1}=C^{-1}BC$ for some 
$C\in G$, and $V^2=[C^{-1}A,B]$, or
\item\label{m=4} $m=4$, $V^2=ABC$ with $A^2=B^2=C^2=1$, and $V^4=[BA,BC]$.
\end{parts}
\end{theorem}

We recall that in a free product, every element of finite order lies in a 
conjugate of a free factor.  Thus we have the following consequence of 
Theorem~\ref{main}.

\begin{corollary}\label{noeven}
Let $G=\freeprod$ where no $G_i$ has elements of even order.  If $V,X,Y\in G$ 
and $V^m=[X,Y]$ for some $m\ge 2$, then either $V\in W^{-1}G_i W$ for some 
$W\in G$, $i\in I$, and $V^m$ is a commutator in $W^{-1}G_i W$ or $m=3$, 
$V=AB$ for some $A,B\in G$ with $A^3=B^3=1$, and $V^3=[BA^2,BAB]$. 
\end{corollary}

Part (\ref{main}.\ref{m=4}) of Theorem \ref{main} is somewhat unsatisfactory 
in that it describes 
the form of $V^2$ rather than that of $V$.  Among the ways in which an element 
$V$ of a free product may have $V^2=ABC$ with $A^2=B^2=C^2=1$ is $V=DE$ 
with $D^2=E^4=1$, in which case $V^2=(D)(E^2)(E^{-1}DE)$.  Not every solution 
is of this form, as shown by 
$G=\langle a,b~;~a^2=b^2=(ab)^2=1\rangle \ast \langle c~;~c^2=1\rangle$ 
and $V=acbcabc$; here $V^2=(acbca)(bcacb)(cabc)$, a product of three elements 
of order two, but $V$ is not a product of two elements of finite order.  A 
classification of elements $V$ satisfying the conditions of 
(\ref{main}.\ref{m=4}) has eluded us.                                  

Relative to (\ref{main}.\ref{m=2}), we record the following well known 
consequence of the Conjugacy Theorem for Free Products 
\cite[Theorem IV.1.4]{LyndonSchupp}:

\begin{lemma}\label{invconj}
If $B$ is an element of a free product $G=\freeprod$ and $B^{-1}=C^{-1}BC$ for 
some $C\in G$, then either
\begin{parts}
\item\label{selfconj} $B\in W^{-1}G_i W$ for some $W\in G$ and $i\in I$ and 
there is a $C\in W^{-1}G_i W$ such that $B^{-1}=C^{-1}BC$, or
\item\label{twoinvol} $B=DE$ for some $D,E\in G$ with $D^2=E^2=1$.
\end{parts}
\end{lemma}

Before proceeding with a proof of the theorem, we establish some notation and 
terminology for the free product $G=\freeprod$.  Our usage is that of 
R.~C.~Lyndon and P.~E.~Schupp \cite{LyndonSchupp} unless otherwise noted.  
A product $PQ$ of elements 
$P$ and $Q$ of $G$ is {\em reduced\/} if one of $P$, $Q$ is trivial or the last 
letter of the normal form of $P$ is not inverse to the first letter of the 
normal form of $Q$.  The product $PQ$ is {\em fully reduced\/} if $P$ or $Q$ 
is trivial or if the last letter of the normal form of $P$ is from a free 
factor different from that of the first letter of the normal form of $Q$; we 
sometimes denote this by writing $P\cdot Q$.  These notions extend to products 
of more than two factors, with the understanding that the non-cancellation 
conditions continue to apply after trivial factors have been deleted.  Thus a 
product $P_1\dots P_k$ is fully reduced if and only if 
$|P_1\dots P_k|=\sum_{i=1}^k |P_i|$, where $|\,\,|$ denotes free product length.

An element $P$ of $G$ is {\em cyclically reduced\/} if $|P|\le 1$ or the first 
and last letters of its normal form are not inverses, and is {\em fully 
cyclically reduced\/} if $|P|\le 1$ or the first and last letters of its 
normal form lie in different free factors of $G$.

A key ingredient in our analysis will be the characterization by M.~J.~Wicks 
of the fully reduced forms of a commutator in a free product.  The following 
is a restatement of Lemma 6 of \cite{Wicks62}.

\begin{lemma}[Wicks]\label{Wicks}
If $U\in G=\freeprod$ is a commutator, either
$U\in W^{-1}G_i W$ for some $W\in G$, $i\in I$, and $U$ 
is a commutator in $W^{-1}G_i W$, or
some fully cyclically reduced conjugate of $U$ has one of the following fully 
reduced forms:
\begin{parts}
\item\label{1var} $X^{-1}a_1Xa_2$ with $X\ne 1$, $a_1\ne 1$, $a_1,a_2\in G_i$ 
for some $i\in I$, and $a_1$ conjugate to $a_2^{-1}$ in $G_i$, or
\item\label{2var} $X^{-1}a_1Y^{-1}a_2Xa_3Ya_4$ with $X\ne 1$, $Y\ne 1$, 
$a_1,a_2,a_3,a_4\in G_i$ for some $i\in I$, and $a_4a_3a_2a_1=1$, or 
\item\label{3var} $X^{-1}a_1Y^{-1}b_1Z^{-1}a_2Xb_2Ya_3Zb_3$ with 
$a_1,a_2,a_3\in G_i$ for some $i\in I$ and $a_3a_2a_1=1$, $b_1,b_2,b_3\in G_j$ 
for some $j\in I$ and $b_3b_2b_1=1$, and either not all of 
$a_1,a_2,a_3,b_1,b_2,b_3$ are in any one free factor of $G$ or each of $X,Y,Z$ 
is nontrivial.
\end{parts}
\end{lemma}

As a final preliminary step, we examine the ways in which both an element and 
its inverse can occur as fully reduced subwords of a proper power in a free 
product.

\begin{lemma}\label{cancel}
Suppose that $V$ is a fully cyclically reduced element of $G=\freeprod$ with 
$|V|\ge 2$, that $m\ge 1$, and that for some $X,R,S,T\in G$, 
$V^m=X^{-1}\cdot R=S\cdot X\cdot T$.  Then one of the following is true:
\begin{parts}
\item\label{verylong} $|X|\ge|V|$, $X=X_1\cdot B\cdot A$ and $V=A\cdot B$ 
for some $A,B,X_1$ with $A^2=B^2=1$, and $SX=V^n\cdot A$ for some $n<m$, 
or
\item\label{long} $\frac{1}{2}|V|<|X|<|V|$, $X=X_1\cdot X_2\cdot X_3$ and
$V=X_3\cdot X_2^{-1}\cdot X_1\cdot X_2$ for some $X_1,X_2,X_3$ with 
$X_1^2=X_3^2=1$, and $S=V^n\cdot X_3\cdot X_2^{-1}$ for some $n<m$, or
\item\label{medium1} $|X|<|V|$, $X=X_1\cdot X_2$ and  
$V=X_2^{-1}\cdot X_1\cdot X_2\cdot T_1$ for some $X_1,X_2,T_1$ 
with $X_1^2=1$, and $S=V^n\cdot X_2^{-1}$ for some $n<m$, or
\item\label{medium2} $|X|<|V|$, $X=X_1\cdot X_2$ and
$V=X_2\cdot X_1^{-1}\cdot S_2\cdot X_1$ for some $X_1,X_2,S_2$ 
with $X_2^2=1$,  
and $S=V^n\cdot X_2\cdot X_1^{-1}\cdot S_3$ for some $n<m$, or
\item\label{short} $|X|\le \frac{1}{2}|V|-1$ and 
$V=X^{-1}\cdot V_2\cdot X\cdot V_3$ for some nontrivial $V_2,V_3$ and 
$S=V^n\cdot X^{-1}\cdot V_2$ for some $n<m$.
\end{parts}
\end{lemma}

\begin{pf*}{Proof of Lemma \ref{cancel}}
If $X$ is empty, clause 
(\ref{cancel}.\ref{short}) applies with $V=V_2\cdot V_3$ a fully reduced 
factorization of $V$ 
such that $S=V^n\cdot V_2$ for some $n<m$.  We suppose, then, that $X\ne 1$.

If $|X|\ge |V|$, we factor $V$ as $A\cdot B$ so that $SX=V^n\cdot A$ with 
$|A|<|V|$.  It follows that $X=X_1\cdot B\cdot A$ for some $X_1$.  But since 
$X^{-1}=A^{-1}\cdot B^{-1}\cdot X_1^{-1}$ is an initial subword of 
$V^m=(A\cdot B)^m$, $A^{-1}=A$ and $B^{-1}=B$.  This is the situation 
described in (\ref{cancel}.\ref{verylong}).  We assume, henceforth, that 
$|X|<|V|$.

Let $n$ be the largest integer such that $|V^n|\le |S|$, and let $S_1$, $V_1$ 
be such that $S=V^n\cdot S_1$ and $V=X^{-1}\cdot V_1$.  We cannot have 
$|S_1|=|X|$ or $|S_1|+|X|=|V|$, for that would violate our hypotheses on the 
fully reduced factorizations of $V^m$.

Suppose that $|S_1|<|X|$ and $|S_1|+|X|>|V|$.  Then $X$ factors as 
$X_1\cdot X_2\cdot X_3$ with $X^{-1}=S_1\cdot X_1^{-1}$, 
$V=S_1\cdot X_1\cdot X_2$, and $X_1$ and $X_2$ nonempty.  Now 
$S_1=X_3^{-1}\cdot X_2^{-1}$, so $V=X_3^{-1}\cdot X_2^{-1}\cdot X_1\cdot X_2$.  
But $SX=V^{n+1}\cdot X_3$, which implies that $X_3^{-1}=X_3$, and 
$V=X_3^{-1}\cdot X_2^{-1}\cdot X_1^{-1}\cdot V_1$, which yields 
$X_1^{-1}=X_1$.  This is the situation of (\ref{cancel}.\ref{long}), and we 
note that $|V|<|S_1|+|X|$ and $|S_1|<|X|$ imply that $|V|<2|X|$.

Next suppose that $|S_1|<|X|$ and $|S_1|+|X|<|V|$.  Then $X$ factors as 
$X_1\cdot X_2$ with $S_1=X_2^{-1}$ and $V=S_1\cdot X\cdot T_1$ for some $T_1$,
and so $V=X_2^{-1}\cdot X_1\cdot X_2\cdot T_1=X_2^{-1}\cdot X_1^{-1}\cdot 
V_1$.  It follows that $X_1^{-1}=X_1$, and we are in situation 
(\ref{cancel}.\ref{medium1}).

Now suppose that $|S_1|>|X|$ and $|S_1|+|X|>|V|$.  We factor $X$ as 
$X_1\cdot X_2$ with $V=S_1\cdot X_1$ and factor $S_1$ as $X^{-1}\cdot S_3$.  
Then $V=X_2^{-1}\cdot X_1^{-1}\cdot S_3\cdot X_1$ and, since 
$S\cdot X=V^{n+1}\cdot X_2$, $X_2^{-1}=X_2$; this is 
(\ref{cancel}.\ref{medium2}).

Finally, suppose that $|S_1|>|X|$ and $|S_1|+|X|<|V|$.  In this case, $S_1$ 
factors as $X^{-1}\cdot V_2$ for some $V_2$ and $V=S_1\cdot X\cdot V_3$ for 
some $V_3$.  Then $V=X^{-1}\cdot V_2\cdot X\cdot V_3$ where necessarily $V_2$ 
and $V_3$ are nonempty, and (\ref{cancel}.\ref{short}) applies.
\end{pf*}

\begin{pf*}{Proof of Theorem \ref{main}}
Each of the forms specified for $V$ (or, in (\ref{main}.\ref{m=4}), $V^2$) in 
the conclusion of Theorem \ref{main} is preserved if $V$ is replaced by a 
conjugate of itself, so we lose no generality in assuming that $V$ is fully 
cyclically reduced.  If $V\in G_i$ for some $i\in I$, then Lemma \ref{Wicks} 
tells us that (\ref{main}.\ref{factor}) holds.  We suppose, then, that 
$|V|\ge 2$.

By Lemma \ref{Wicks}, some fully cyclically reduced conjugate of $V^m$ has the 
form specified in (\ref{Wicks}.\ref{1var}), (\ref{Wicks}.\ref{2var}), or
(\ref{Wicks}.\ref{3var}).  After again replacing $V$ by a fully cyclically 
reduced conjugate and relabeling in (\ref{Wicks}.\ref{2var}) and 
(\ref{Wicks}.\ref{3var}) if necessary, we may assume that $V^m$ has form 
(\ref{Wicks}.\ref{1var}), or form (\ref{Wicks}.\ref{2var}) with 
$|X|\ge |Y|$, or form (\ref{Wicks}.\ref{3var}) with $|X|\ge |Y|$ and 
$|X|\ge |Z|$.  

Let $P=a_1$ and $Q=a_2$ in form (\ref{Wicks}.\ref{1var}), $P=a_1Y^{-1}a_2$ and 
$Q=a_3Ya_4=\allowbreak
a_3Ya_1^{-1}a_2^{-1}a_3^{-1}$ in form (\ref{Wicks}.\ref{2var}), and 
$P=a_1Y^{-1}b_1Z^{-1}a_2$ and $Q=b_2Ya_3Zb_3=\allowbreak
b_2Ya_1^{-1}a_2^{-1}Zb_1^{-1}b_2^{-1}$ in 
form (\ref{Wicks}.\ref{3var}). In each 
instance, $V^m=X^{-1}\cdot P\cdot X\cdot Q$ and $Q$ is conjugate to $P^{-1}$ in 
$G$.  Further, $|P|=|Q|=1$ in (\ref{Wicks}.\ref{1var}), $|P|\le |X|+2$ and 
$|Q|\le |X|+2$ in (\ref{Wicks}.\ref{2var}), and $|P|\le 2|X|+3$ and 
$|Q|\le 2|X|+3$ in (\ref{Wicks}.\ref{3var}).  We proceed by cases according to 
which clause of the conclusion of Lemma \ref{cancel} is satisfied, with 
$R=PXQ$, $S=X^{-1}P$, and $T=Q$.

\subsection*{Case \rom{(\ref{cancel}.\ref{verylong})}}
Suppose that $X=X_1\cdot B\cdot A$ and $V=A\cdot B$ for some $X_1,A,B$ with 
$A^2=B^2=1$, that $X_1^{-1}PX_1=(AB)^kA$ for some $k$, $0\le k\le m-3$, and 
that $Q=B(AB)^{m-k-3}$.

If $m$ is even, (\ref{main}.\ref{m=even}) is satisfied, while if $m$ is odd, 
$Q$ conjugate to $P^{-1}$ implies that $B$ is conjugate to $A$ and 
(\ref{main}.\ref{m=odd}) holds.

\subsection*{Case \rom{(\ref{cancel}.\ref{long})}}
Suppose that $X=X_1\cdot X_2\cdot X_3$ and $V=X_3\cdot X_2^{-1}\cdot X_1\cdot 
X_2$ for some $X_1,X_2,X_3$ with $X_1^2=X_3^2=1$, that 
$P=X_2X_3X_2^{-1}(X_1X_2X_3X_2^{-1})^k$ for some $k$, $0\le k\le m-3$, and 
$Q=X_2^{-1}X_1X_2(X_3X_2^{-1}X_1X_2)^{m-k-3}$.

As in the previous case, (\ref{main}.\ref{m=even}) applies if $m$ is even, and 
if $m$ is odd, $Q$ conjugate to $P^{-1}$ implies that $X_3$ is conjugate to 
$X_1$ and so (\ref{main}.\ref{m=odd}) obtains.

\subsection*{Case \rom{(\ref{cancel}.\ref{medium1})}}
Suppose that $|X|<|V|$, $X=X_1\cdot X_2$ and $V=X_2^{-1}\cdot X_1\cdot 
X_2\cdot T_1$ for some $X_1,X_2,T_1$ with $X_1^2=1$, that 
$P=X_2T_1X_2^{-1}(X_1X_2T_1X_2^{-1})^k$ for some $k$, $0\le k\le m-2$, and 
$Q=T_1(X_2^{-1}X_1X_2T_1)^{m-k-2}$.  

We first notice that since $|P|\le 2|X|+3\le 2|V|+1$ and 
$|Q|\le 2|X|+3\le 2|V|+1$, we have $m\le 6$.  Now $Q$ is conjugate to 
$P^{-1}$, so $P$ and $Q$ must have fully cyclically reduced conjugates of the 
same length.  It is not hard to see that this implies that either $k=m-k-2$ or 
$T_1^2=1$.  If $T_1^2=1$, we find as in previous cases that 
(\ref{main}.\ref{m=even}) applies if $m$ is even and that 
(\ref{main}.\ref{m=odd}) applies if $m$ is odd.  We suppose, then, that 
$T_1^2\ne 1$ and that $k=m-k-2$.  The possibilities to consider are that $m=2$ 
and $k=0$, $m=4$ and $k=1$, and $m=6$ and $k=2$.

If $m=2$ and $k=0$, $T_1$ is conjugate to $T_1^{-1}$ and 
(\ref{main}.\ref{m=2}) holds.

If $m=4$ and $k=1$, $Q=T_1X_2^{-1}X_1X_2T_1$ and 
$P=X_2T_1X_2^{-1}X_1X_2T_1X_2^{-1}$, a conjugate of $Q$.  Now $T_1^2\ne 1$, so 
$Q$ is not in a conjugate of a free factor of $G$, but since $Q$ is conjugate 
to $P^{-1}$, $Q$ is conjugate to $Q^{-1}$.  By Lemma \ref{invconj}, then, 
$Q=DE$ for some $D,E$ with $D^2=E^2=1$.  But then 
$V^2=X_2^{-1}X_1X_2DE$ and (\ref{main}.\ref{m=4}) applies.

Suppose, then, that $m=6$ and $k=2$.  We must have $|X|=|V|-1$ and 
$|P|=|Q|=2|V|+1$, so $X_2$ is empty and $T_1$ has length one.  Let us write 
$X_1=C^{-1}\cdot a\cdot C$ with $C\in G$ and $a\in G_i$ for some $i\in I$ and 
$a^2=1$ and $T_1=b\in G_j$ for some $j\in I$ with $b^2\ne 1$.  We then have 
$P=Q=b\cdot C^{-1}\cdot a\cdot C\cdot b\cdot C^{-1}\cdot a\cdot C\cdot b$, and 
so $b^2\cdot C^{-1}\cdot a\cdot C\cdot b\cdot C^{-1}\cdot a\cdot C$ is a fully 
cyclically reduced conjugate of $P$ which,
like $P$, is conjugate to its inverse.  There 
must then be a factorization $C_1\cdot C_2$ of $C$ such that one of the 
following holds:
\begin{equation}
C_1^{-1}aC_1C_2b^{-1}C_2^{-1}C_1^{-1}aC_1C_2b^{-2}C_2^{-1}=
b^2C_2^{-1}C_1^{-1}aC_1C_2bC_2^{-1}C_1^{-1}aC_1C_2, \label{one} 
\end{equation}
\begin{equation}
C_2b^{-1}C_2^{-1}C_1^{-1}aC_1C_2b^{-2}C_2^{-1}C_1^{-1}aC_1=
b^2C_2^{-1}C_1^{-1}aC_1C_2bC_2^{-1}C_1^{-1}aC_1C_2, \label{two}     
\end{equation}
\begin{equation}
C_1^{-1}aC_1C_2b^{-2}C_2^{-2}C_1^{-1}aC_1C_2b^{-1}C_2^{-1}=
b^2C_2^{-1}C_1^{-1}aC_1C_2bC_2^{-1}C_1^{-1}aC_1C_2, \label{three} 
\end{equation}
\begin{equation} 
C_2b^{-2}C_2^{-1}C_1^{-1}aC_1C_2b^{-1}C_2^{-1}C_1^{-1}aC_1=
b^2C_2^{-1}C_1^{-1}aC_1C_2bC_2^{-1}C_1^{-1}aC_1C_2. \label{four} 
\end{equation}

If (\ref{one}) is true, a length comparison on the fully reduced products on 
the two sides shows that
\[ C_1^{-1}aC_1C_2b^{-1}C_2^{-1}=b^2C_2^{-1}C_1^{-1}aC_1C_2\]
and
\[ C_1^{-1}aC_1C_2b^{-2}C_2^{-1}=bC_2^{-1}C_1^{-1}aC_1C_2.\]
The left sides of these two equations begin with the same normal form letter, 
so looking at the right sides we get $b^2=b$, a contradiction.  Similarly, 
(\ref{two}) yields
\[ C_2b^{-1}C_2^{-1}C_1^{-1}aC_1=b^2C_2^{-1}C_1^{-1}aC_1C_2\]
and 
\[ C_2b^{-2}C_2^{-1}C_1^{-1}aC_1=bC_2^{-1}C_1^{-1}aC_1C_2,\]
from which we get the contradiction $b^2=b$ if $C_2$ is nonempty, or the 
equation $b^{-1}=b^2$ if $C_2$ is empty.  This last possibility corresponds 
to (\ref{main}.\ref{m=6}).  If (\ref{three}) holds, we get
\[ C_1^{-1}aC_1C_2b^{-2}C_2^{-1}=b^2C_2^{-1}C_1^{-1}aC_1C_2\]
and
\[ C_1^{-1}aC_1C_2b^{-1}C_2^{-1}=bC_2^{-1}C_1^{-1}aC_1C_2.\]
As in (\ref{one}), we derive the contradiction $b^2=b$.  Finally, if 
(\ref{four}) is true, 
\[ C_2b^{-2}C_2^{-1}C_1^{-1}aC_1=b^2C_2^{-1}C_1^{-1}aC_1C_2\]
and
\[ C_2b^{-1}C_2^{-1}C_1^{-1}aC_1=bC_2^{-1}C_1^{-1}aC_1C_2.\]
This yields the contradictions $b^{-1}=b$ if $C_2$ is empty and $b^2=b$ if 
$C_2$ is nonempty.

\subsection*{Case \rom{(\ref{cancel}.\ref{medium2})}}
Suppose that $|X|<|V|$, $X=X_1\cdot X_2$ and $V=X_2\cdot X_1^{-1}\cdot 
S_2\cdot X_1$ for some $X_1,X_2,S_2$ with $X_2^2=1$, that 
$P=S_2(X_1X_2X_1^{-1}S_2)^k$ for some $k$, $0\le k\le m-2$, and 
$Q=X_1^{-1}S_2X_1(X_2X_1^{-1}S_2X_1)^{m-k-2}$.

Replacing $V$ by its fully cyclically reduced conjugate
$X_1X_2X_1^{-1}S_2$ and changing 
notation reduces this to Case (\ref{cancel}.\ref{medium1}).

\subsection*{Case \rom{(\ref{cancel}.\ref{short})}}
Suppose that $|X|\le \frac{1}{2}|V|-1$, $V=X^{-1}\cdot V_2\cdot X\cdot V_3$ 
for some $V_2,V_3$, that $P=V_2(XV_3X^{-1}V_2)^k$ for some $k$, 
$0\le k\le m-1$, and that $Q=V_3(X^{-1}V_2XV_3)^{m-k-1}$.

Since $|P|\le 2|X|+3\le |V|+1$ and $|Q|\le 2|X|+3\le |V|+1$, we have $m\le 3$.  
We first consider the case that $m=2$.  If $k=0$, $Q=V_3X^{-1}V_2XV_3$ 
conjugate to $P^{-1}=V_2^{-1}$ implies that $V_3^2=1$ and $V_2$ is conjugate 
to $V_2^{-1}$; (\ref{main}.\ref{m=2}) applies.  If $k=1$, 
$P=V_2XV_3X^{-1}V_2$ is conjugate to $Q^{-1}=V_3^{-1}$, so $V_2^2=1$, $V_3$ is 
conjugate to $V_3^{-1}$, and again (\ref{main}.\ref{m=2}) applies.

Now suppose that $m=3$.  In this event, we must have $|X|=\frac{1}{2}|V|-1$ 
and $|P|=|Q|=|V|+1$, so $|V_2|=|V_3|=1$.  Let us write $V_2=a\in G_i$ for some 
$i\in I$ and $V_3=b\in G_j$ for some $j\in I$.  Then since 
$Q=b(X^{-1}aXb)^{2-k}$ is conjugate to $P^{-1}=a^{-1}(Xb^{-1}X^{-1}a^{-1})^k$, 
either $a^2=b^2=1$ and $a$ is 
conjugate to $b$, as described in (\ref{main}.\ref{m=odd}), or $a^2\ne 1$, 
$b^2\ne 1$, $k=1$, and there is a factorization $X_1\cdot X_2$ of $X$ such that 
one of the following holds:
\begin{equation}
X_2b^{-1}X_2^{-1}X_1^{-1}a^{-2}X_1=
b^2X_2^{-1}X_1^{-1}aX_1X_2,\label{first} 
\end{equation}
\begin{equation}
X_1^{-1}a^{-2}X_1X_2b^{-1}X_2^{-1}=
b^2X_2^{-1}X_1^{-1}aX_1X_2.\label{second}
\end{equation}
If (\ref{first}) is true, either $X_2$ is empty and $a^3=b^3=1$ as in 
(\ref{main}.\ref{m=3}), or $X_2$ is nonempty and $X_2b^{-1}=b^2X_2^{-1}$, so 
that $X_2=b^2X_3$ and $X_2^{-1}=X_3^{-1}b^{-1}$ for some $X_3$, producing the 
contradiction $b^2=b$.  If (\ref{second}) is true, $X_2^2=1$ and 
\[ X_1^{-1}a^{-2}X_1X_2b^{-1}=b^2X_2^{-1}X_1^{-1}aX_1.\]
If $X_1$ is nonempty, $X_1=X_4b^{-1}$ and $X_1^{-1}=b^2X_4^{-1}$ for some 
$X_4$, whence $b^{-1}=b^{-2}$, a contradiction.  Thus $X_1$ is empty, and 
$a^{-2}X_2b^{-1}=b^2X_2^{-1}a$ implies that $b^{-1}=a$ and $X_2=X_2^{-1}$.  
Thus $V=XaXa^{-1}$ with $X^2=1$, and (\ref{main}.\ref{m=odd}) applies.
\end{pf*}

\end{document}